\documentclass[12pt]{article}

\usepackage{graphicx}        
\usepackage[bottom]{footmisc}
\usepackage{algorithm}
\usepackage{algpseudocode}

\usepackage{amssymb} 
\usepackage{amsfonts} 

\usepackage{url}

\newtheorem{defn}{Definition}[section]

\newtheorem{example}[defn]{Example}
\newtheorem{theorem}[defn]{Theorem}
\newtheorem{proposition}[defn]{Proposition}
\newtheorem{corollary}[defn]{Corollary}

\def\bp{\noindent{\bf Proof. \ }}
\def\ep{\noindent{$\Box$}}

\begin{document}

\title{A survey on the unconditional convergence \\ and the invertibility of multipliers \\ with implementation}

\author{\vspace{.1in} Diana T. Stoeva and Peter Balazs \\
{}
 \vspace{-.07in}  
 {\small   Acoustics Research Institute, Austrian Academy of Sciences,}
\\
\vspace{-.07in} 
{\small  Wohllebengasse 12-14, Vienna 1040, Austria}\\
{\small  Emails: dstoeva@kfs.oeaw.ac.at, peter.balazs@oeaw.ac.at
}
}

\maketitle

\abstract{
The paper presents a survey over frame multipliers and related concepts. In particular, it includes a short motivation of why multipliers are of interest to consider, a review as well as extension of recent results, devoted to the unconditional convergence of multipliers, sufficient and/or necessary conditions for the invertibility of multipliers, and representation of the inverse via Neumann-like series and via multipliers with particular parameters. Multipliers for frames with specific structure, namely, Gabor and wavelet multipliers, are also considered. Some of the results for the representation of the inverse multiplier are implemented in Matlab codes and the algorithms are described.}

\vspace{.1in} \noindent
{\bf Keywords} Multiplier; Gabor multiplier; Frame; Dual frame; Invertibility; Unconditional convergence.

\section{Introduction}\label{Sec:introd}

Multipliers are operators which consist of an analysis stage, a multiplication, and then a synthesis stage (see Definition \ref{definmult}). 
This is a very natural concept, that occurs in a lot of scientific questions in mathematics, physics, and engineering.

In {\em Physics}, multipliers are the link between classical and quantum mechanics, so called quantization operators \cite{aliant1}. Here multipliers link sequences (or functions) $m_k$ corresponding to the measurable variables in classical physics, to operators $M_{m,\Phi,\Psi}$, which are the measurables in quantum mechanics, via (\ref{defmult}), see e.g. 
\cite{Gaz09,colgaz10}. 

In {\em Signal Processing}, multipliers are a particular way to implement time-variant filters \cite{hlawatgabfilt1}. 
One of the goals in signal processing is to find a transform which
allows certain properties of the signal to be easily found or seen. 
Via such a transform, one can focus on those properties of the signal, one is interested in, or would like to change.
The coefficients can be manipulated directly in the transform domain and thus, certain signal features can be amplified or attenuated. This is for example the case of what a sound engineer does during a concert, operating an equalizer, i.e. changing the amplification of certain frequency bands in real time. 
Filters, i.e. convolution operators, correspond to a multiplication in the Fourier domain, and therefore to a time-invariant change of frequency content. 
They are one of the most important concepts in signal processing.  
Many approaches in signal processing assume a quasi-stationary assumption, i.e., a shift-invariant approach is assumed only locally. There are several ways to have a true time-variant approach (while still keeping the relation to the filtering concept), and one of them is to use the so called Gabor Filters \cite{xxlframehs07,hlawatgabfilt1}. These are Gabor multipliers, i.e. time-frequency multipliers. 

An additional audio signal processing application is the transformation of a melody played by one instrument to sound like played by another.
For an approach of how to identify multipliers which transfer one given signal to another one, see \cite{oltokr13}.

In {\em Acoustics}, the time-frequency filters are used in several fields, for example in Computational Auditory Scene Analysis (CASA) \cite{wanbro06}. The CASA-term refers to the study of auditory scene analysis by computational means, i.e. the separation of auditory events. The CASA problem is closely related to the problem of source separation. Typically, an auditory-based time frequency transform  is calculated from certain  acoustic features and so-called ``time-frequency masks'' are generated. These masks are directly applied to the perceptual representation; they weight the ``target'' regions (mask = 1) and suppress the background (mask = 0). This corresponds to a binary time-frequency multipliers. Such adaptive filters are also used in perceptual sparsity, where a time-frequency mask is estimated from the signal and a simple psychoacustical masking model, resulting in a time-frequency transform, where perceptual irrelevant coefficients are deleted, see \cite{xxllabmask1,necc11}.

Last but certainly not least, multipliers were and are of utmost importance in {\em Mathematics}, where they are used for the diagonalization of operators. Schatten used multipliers for orthonormal bases to describe certain classes of operators \cite{schatt1}, later known as Schatten $p$-classes. The well-known spectral theorems, see e.g. \cite{conw1}, are just results stating that certain operators can always be represented as multipliers of orthonormal bases. Because of their importance for signal processing, Gabor multipliers were defined as time-frequency multipliers \cite{feinow1,benepfand07,groech10}, which motivated the definition of multipliers for general frames in \cite{xxlmult1}.
Recently, the formal definition of frame multipliers led to a  
lot of new approaches to multipliers  \cite{rahxxlXX,Arias2008581,rahgenmul10,xxlbayasg11} and new results \cite{balsto09new,bstable09,Futamura20123201}. 

\begin{figure}[ht] 
\begin{center}
\includegraphics[width=1\textwidth]{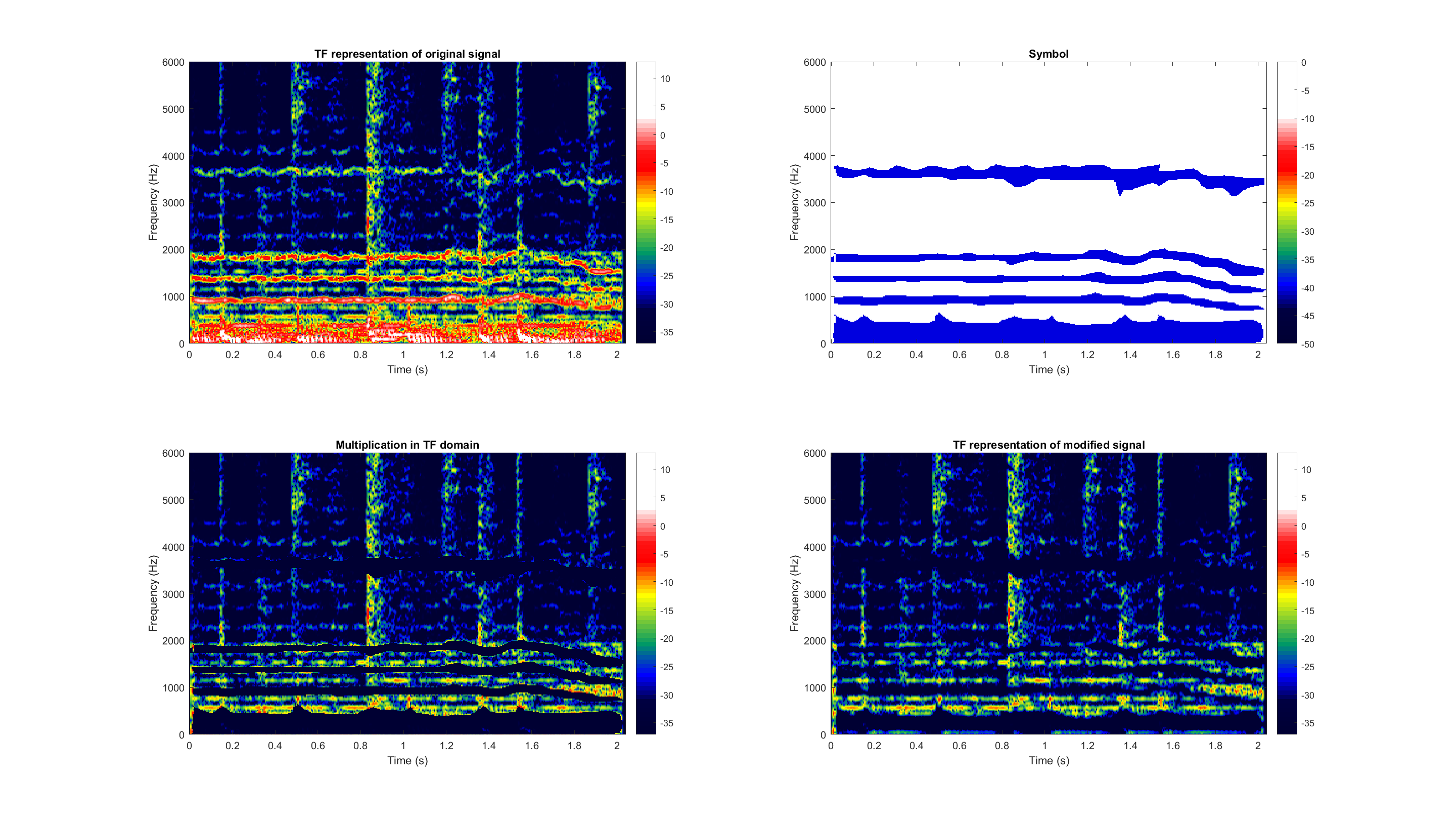} 
\caption{\label{fig:examp1}\small {\em An illustrative example to visualize a multiplier.} 
(TOP LEFT) The time-frequency representation of 
 the music signal $f$. (TOP RIGHT) The symbol $m$, found by a (manual) estimation of the time-frequency region of the singer's voice. 
 (BOTTOM LEFT) The multiplication in the TF domain.
(BOTTOM RIGHT) Time-frequency representation of $M_{m,\widetilde \Psi,\Psi}f$.} 
\end{center}
\end{figure}

Like in \cite{bsreprinv2015}  we show a visualization of a multiplier $M_{m,\widetilde{\Psi},\Psi}$ in the time-frequency plane in Figure \ref{fig:examp1}, using a different setting at a different soundfile.
The visualization is done using algorithms in the LTFAT toolbox \cite{soendxxl10}, in particular using the graphical user interface {\em MULACLAB} \footnote{\url{http://ltfat.github.io/doc/base/mulaclab.html}} for multipliers \cite{ltfatnote030}.  We consider $2$ seconds long excerpt of a Bulgarian folklore song ("Prituri se planinata", performed by Stefka Sabotinova) as signal $f$.
For a time-frequency representation of the musical signal $f$ (TOP LEFT) we use a Gabor frame $\Psi$ (a $23$ ms Hanning window with $75$\% overlap and double length FFT). This $\Psi$ constitutes a Parseval Gabor frame, so it is self-dual. 
By a manual estimation, we determine the symbol $m$ that should describe the time-frequency region of the singer's voice. This region is then multiplied by $0.01$, the rest by $1$ (TOP RIGHT). 
Finally, we show the multiplication in the TF domain (BOTTOM LEFT), and time-frequency representation of the modified signal (BOTTOM RIGHT). 

In this paper we deal with the mathematical concept of multipliers, in particular with frame multipliers.
We will give a survey over the mathematical properties of these operators, collecting known results and combining them with new findings and giving accompanying implementation.
We consider the case of multipliers for general sequences (Section \ref{Sec:2generalmultipl}), as well as multipliers for Gabor and wavelet systems (Section \ref{Sec:exampl}). The paper is organized as follows.  In Section \ref{secNotation} we collect the  basic needed definitions and the notation used in the paper. 
 In Section \ref{Sec:2generalmultipl}, first we discuss well-definedness of multipliers, as well as necessary and sufficient conditions for the unconditional convergence of certain classes of multipliers. Furthermore, we list relations between the symbol of a multiplier and the operator-type of the multiplier. Finally, we consider injectivity, surjectivity, and invertibility of multipliers, presenting necessary and/or sufficient conditions. 
In the cases of invertibility of frame multipliers, we give formulas for the inverse operators in two ways - as Neumann-like series and as multipliers determined by the reciprocal symbol and appropriate dual frames of the initial ones. For the formulas of the inverses given as Neumann-like series, we provide implementations via Matlab-codes. In Section \ref{Sec:exampl}, first we state consequences of the general results on unconditional convergence applied to Gabor and wavelet systems; next we consider invertibility of Gabor multipliers and representation of the inverses via Gabor multipliers with dual Gabor frames of the initial ones.
Finally, in Section \ref{sec-implement}, we implement the inversion of frame multipliers according to some of the statements in Section \ref{Sec:2generalmultipl} and visualize the convergence rate of one of the algorithms in Fig. \ref{fig:ConvRateGab}. 
For the codes of the implementations, as well as for the scripts and the source files which were used to create Fig. \ref{fig:examp1} and Fig. \ref{fig:ConvRateGab}, see the webpage
\url{https://www.kfs.oeaw.ac.at/InversionOfFrameMultipliers}.

\section{Notation and basic definitions}\label{secNotation}

Throughout the paper, $\mathcal{H}$ denotes a separable Hilbert space and $I$ denotes a countable index set. If not stated otherwise, $\Phi$ and $\Psi$ denote sequences $(\phi_n)_{n\in I}$ and $(\psi_n)_{n\in I}$, respectively, having elements from $\mathcal{H}$; $m$ denotes a complex number sequence $(m_n)_{n\in I}$, $\overline{m}$ - the sequence of the complex conjugates of $m_n$, $1/m$ - the sequence of the reciprocals on $m_n$, $m\Phi$ - the sequence $(m_n \phi_n)_{n\in I}$. 
The sequence $\Phi$ (resp. $m$) is called {\em norm-bounded below} (in short, $NBB$) if $(\|\phi_n\|_{n\in I})\in\ell^\infty$ (resp. $(|m_n|_{n\in I})\in\ell^\infty$).
An operator $F:\mathcal{H}\to\mathcal{H}$ is called {\it invertible on $\mathcal{H}$} (or just {\it invertible}, if there is no risk of confusion of the space) if it is a bounded bijective operator from $\mathcal{H}$ onto 
$\mathcal{H}$.

Recall that $\Phi$ is called

- a {\it Bessel sequence in $\mathcal{H}$} if there is $B_\Phi\in(0,\infty)$ (called a {\it Bessel bound of $\Phi$}) so that $\sum_{n\in I} |\langle f, \phi_n\rangle |^2\leq B_\Phi \|f\|^2$ for every $f\in\mathcal{H}$;

- a {\it frame for $\mathcal{H}$} \cite{ds52} if there exist  $A_\Phi\in (0,\infty)$ and $B_\Phi\in (0,\infty)$ (called {\it frame bounds of $\Phi$}) so that $A_\Phi \|f\|^2\leq \sum_{n\in I} |\langle f, \phi_n\rangle |^2\leq B_\Phi \|f\|^2$ for every $f\in\mathcal{H}$;

- a {\it Riesz basis for $\mathcal{H}$} \cite{bari51} if it is the image of an orthonormal basis under an invertible operator.

Recall the needed basics from frame theory (see e.g. \cite{ole1}).
Let $\Phi$ be a frame for $\mathcal{H}$. Then there exists a frame $\Psi$ for $\mathcal{H}$ so that $f=\sum_{n\in I} \langle f, \phi_n\rangle \psi_n
= \sum_{n\in I} \langle f, \psi_n\rangle \phi_n$ for every $f\in\mathcal{H}$; such a frame $\Psi$ is called a {\it dual frame of $\Phi$}. 
One associates the analysis operator  $U_\Phi: \mathcal{H}\to\ell^2$ given by $U_\Phi f = (\langle f, \phi_n\rangle)$, the synthesis operator $T_\Phi:\ell^2 \to \mathcal{H}$ given by $T_\Phi (c_n)_{n\in I} = \sum_{n\in I} c_n \phi_n$, and the frame operator $S_\Phi:\mathcal{H}\to \mathcal{H}$ given by $S_\Phi f=\sum_{n\in I} \langle f, \phi_n\rangle \phi_n$; all these operators are well defined and bounded. The frame operator $S_\Phi$ is invertible on $\mathcal{H}$ and the sequence $(S_\Phi^{-1} \phi_n)_{n\in I}$  forms a dual frame of $\Phi$ called {\it the canonical dual}. If $\Phi$ is a Riesz basis for $\mathcal{H}$ (and thus a frame for $\mathcal{H}$), the canonical dual is the only dual frame of $\Phi$. If $\Phi$ is a frame for $\mathcal{H}$ which is not a Riesz basis for $\mathcal{H}$ (so called {\it overcomplete frame}), then in addition to the canonical dual there are other dual frames.  If a sequence $\Psi$ not necessarily being a frame for $\mathcal{H}$  satisfies 
 $f=\sum_{n\in I} \langle f, \phi_n\rangle \psi_n$
  (resp.  $f=\sum_{n\in I} \langle f, \psi_n\rangle \phi_n$)
 for all $f\in\mathcal{H}$, it is called a {\it synthesis} (resp. {\it analysis})
 {\it pseudo-dual of $\Phi$}; for more on analysis and synthesis pseudo-duals see \cite{Scharactarxiv}. 
When a frame $\Psi$ is not necessarily a dual one of $\Phi$, 
but there is $\varepsilon\in(0,1)$ so that $\|f-\sum_{n\in I} \langle f,\phi_n\rangle\psi_n\| \leq \varepsilon \|f\|$ for every $f\in\mathcal{H}$, then 
$\Psi$ is called an {\it approximately dual frame of $\Phi$} \cite{chla08}.

  Let $\Lambda=\{ (\omega, \tau)\}$  be a {\it lattice} in $\mathbb{R}^{2d}$, i.e., a discrete subgroup of $\mathbb{R}^{2d}$ of the form $A\mathbb{Z}^{2d}$ for some invertible matrix $A$. 
For $\omega \in \mathbb{R}^d$ and $\tau\in\mathbb{R}^d$, recall the modulation operator $E_\omega:L^2(\mathbb{R}^d) \to L^2(\mathbb{R}^d)$ determined by 
$ (E_\omega f) (x)=e^{2\pi {\rm i} \omega x} f(x)$ 
and the translation operator $T_\tau:L^2(\mathbb{R}^d)\to L^2(\mathbb{R}^d)$ given by
$(T_\tau f) (x)= f(x-\tau)$.  
For $g\in L^2(\mathbb{R}^d)$, a system (resp. frame  for $L^2(\mathbb{R}^d)$) of the form $(E_\omega T_\tau g)_{(\omega,\tau)\in\Lambda} $ is called a {\emph Gabor system} (resp. \emph{Gabor frame} \emph{for $L^2(\mathbb{R}^d)$}). 
Recall also the dilation operator $D_a:L^2(\mathbb{R}) \to L^2(\mathbb{R})$ (for $a\neq 0$) determined by 
$ (D_a f) (x)=\frac{1}{\sqrt{|a|}} f(\frac{x}{a})$. 
Given $\psi\in L^2(\mathbb{R}^d)$, a system (resp. frame  for $L^2(\mathbb{R}^d)$) of the form $(T_b D_a \psi)_{(a,b)\in\Lambda} $ is called a {\emph wavelet system} (resp. \emph{wavelet frame} \emph{for $L^2(\mathbb{R}^d)$}).

\begin{defn} \label{definmult}
Given $m$, $\Phi$, and $\Psi$, the operator $M_{m,\Phi,\Psi}$ given by
\begin{equation}\label{defmult}
M_{m,\Phi, \Psi} f = \sum m_n \langle f, \psi_n\rangle \phi_n, \ f\in {\mathcal D}(M_{m,\Phi, \Psi}),
\end{equation}
where  ${\mathcal D}(M_{m,\Phi, \Psi})=\{f\in\mathcal{H} \ : \   \sum m_n \langle f, \psi_n\rangle \phi_n \mbox{ converges}\}$, 
is called a \emph{multiplier}. The sequence $m$ is called the \emph{symbol} (also, the \emph{weight}) of the multiplier. When $\Phi$ and $\Psi$ are Gabor systems (resp. Bessel sequences, frames, Riesz bases), $M_{m,\Phi, \Psi}$ is called a \emph{Gabor} (resp. \emph{Bessel, frame, Riesz}) \emph{multiplier}. 
A multiplier $M_{m,\Phi, \Psi}$ is called  \emph{well-defined} (resp. \emph{unconditionally convergent}) if the series in 
(\ref{defmult}) is convergent (resp. unconditionally convergent) in $\mathcal{H}$ for every $f\in \mathcal{H}$. 
\end{defn}

\section{Properties of multipliers}\label{Sec:2generalmultipl} 

In this section we consider  some mathematical properties of multipliers using sequences, 
without assuming a special structure for them, while Section \ref{Sec:exampl} is devoted to certain classes of structured sequences.  

\subsection{On well-definedness, unconditional convergence, and relation to the symbol}
\label{subs-UncondConv}

Multipliers are motivated by applications and they can be intuitively understood.
Still the general definition given alone is a mathematical one, and so we have to clarify basic mathematical properties like boundedness and unconditional convergence.

We collect and state some basic results: 

\begin{proposition}\label{propgen1} For a general multiplier, the following relations to unconditional convergence hold.
  
\begin{itemize}
\item[{\rm (a)}] 
A well-defined multiplier is always bounded, but not necessarily unconditionally convergent.

\item[{\rm (b)}]  
  If $m\in\ell^\infty$,  
then a Bessel multiplier $M_{m,\Phi,\Psi}$ is well-defined 
bounded operator with $\|M_{m,\Phi,\Psi}\|\leq \sqrt{B_\Phi B_\Psi}\, \|m\|_\infty$ and the series in (\ref{defmult}) converges unconditionally for every $f\in \mathcal{H}$. The converse does not hold in general, even for a frame multiplier $M_{m,\Phi,\Psi}$.

\item[{\rm (c)}]
  If $\Phi$ and $\Psi$ are NBB, then a Bessel multiplier $M_{m,\Phi,\Psi}$ is unconditionally convergent if and only if $m\in\ell^\infty$.

\item[{\rm (d)}]
If $\Phi$, $\Psi$, and $m$ are NBB, then  $M_{m,\Phi,\Psi}$ is unconditionally convergent if and only if 
 $M_{m,\Phi,\Psi}$ is a Bessel multiplier and $m$ is semi-normalized.

  \item[{\rm (e)}] 
  If $\Phi$ is NBB Bessel, then $M_{m,\Phi,\Psi}$ is unconditionally convergent if and only if
$m\Psi$ is Bessel.
  \item[{\rm (f)}]
  If $\Phi$ is NBB Bessel and $m$ is NBB, then $M_{m,\Phi,\Psi}$ is unconditionally convergent if and only if 
$M_{m,\Phi,\Psi}$ is a Bessel multiplier.

  \item[{\rm (g)}]
  If $\Phi$ is a Riesz basis for $\mathcal{H}$, then $M_{m,\Phi,\Psi}$ is unconditionally convergent if and only if
it is well defined if and only if $m\Psi$ is Bessel.

  \item[{\rm (h)}]
  If $\Phi$ is a Riesz basis for $\mathcal{H}$ and $\Psi$ is NBB, then well-definedness of $M_{m,\Phi,\Psi}$ implies $m\in\ell^\infty$, 
	while the converse does not hold in general.

\item[{\rm (k)}] 
A Riesz multiplier $M_{m,\Phi,\Psi}$ is well-defined if and only if $m\in\ell^\infty$ if and only if it is unconditionally convergent .

\end{itemize}
\end{proposition}
\bp 
(a) That every well defined multiplier is bounded is stated in \cite[Lemma 2.3]{uncconv2011}, the result follows simply by the Uniform Boundedness Principle. However, not every well defined multiplier is unconditionally convergent - consider for example the sequences in \cite[Remark 4.10(a)]{bstable09}, namely, $\Phi=(e_1, e_2, e_2, e_2, e_3, e_3, e_3, e_3, e_3, \ldots)$ and $\Psi=(e_1, e_2, e_2, -e_2, e_3, e_3, -e_3, e_3, -e_3,\ldots)$ for which one has that $M_{(1),\Phi,\Psi}$ is the identity operator on $\mathcal{H}$ but it is not unconditionally convergent.

\sloppy
(b) The first part of (b) is given in \cite{xxlmult1}.  To show that the converse is not valid in general, i.e., that unconditionally convergent frame multiplier does not require $m\in\ell^\infty$, consider for example \cite[Ex. 4.6.3(iv)]{bstable09}, namely,  $m=(1,1,1,1,2,2,1,3,3,\ldots)$, $\Phi=(e_1, e_1, -e_1, e_2, e_2, -e_2, e_3, e_3, -e_3,\ldots)$, and $\Psi=(e_1, e_1, e_1, e_2, \frac{1}{2} e_2, \frac{1}{2} e_2,e_3, \frac{1}{3} e_3, \frac{1}{3} e_3,\ldots)$, for which one has that $M_{m,\Phi,\Psi}$ is unconditionally convergent.

(c) (resp. (d)) One of the directions follows from (b) and the other one from \cite[Prop. 3.2(iii)]{uncconv2011} (resp. \cite[Prop.3.2(iv)]{uncconv2011}).

(e) is given in \cite[Prop. 3.4(i)]{uncconv2011}.

(f) Let $\Phi$ be NBB Bessel and let $m$ be NBB. If $M_{m,\Phi,\Psi}$ is unconditionally convergent, then by (e) the sequence $m\Psi$ is Bessel which by the NBB-property of $m$ implies that $\Phi$ is Bessel. The converse statement follows from (b). 

(g) and (h) can be found in \cite[Prop. 3.4]{uncconv2011}. 

(k) The first equivalence is given in \cite[Prop. 3.4(iv)]{uncconv2011} and for the second equivalence see (g). 
\ep

A natural question for any linear operator is how its adjoint looks. 
For multipliers we can show the following.

 \begin{proposition}  \label{lemuncb} {\rm \cite{uncconv2011} } 
 For any $\Phi,\Psi$ and $m$, the following holds. 
\begin{itemize}
\item[{\rm (i)}] 
If $M_{m,\Phi,\Psi}$ is well defined (and hence bounded), then  $M_{m,\Phi,\Psi}^*$ equals $M_{\overline{m},\Psi,\Phi}$  in a weak sense. 
\item[{\rm (ii)}] If $M_{m,\Phi,\Psi}$ and  $M_{\overline{m},\Psi,\Phi}$ are well defined on all of $\mathcal{H}$, then $M_{m,\Phi,\Psi}^*=M_{\overline{m},\Psi,\Phi}$.
\end{itemize}
\end{proposition}

For more on well definedness and unconditional convergence, we refer to \cite{uncconv2011}.  \\

As a consequence of Proposition \ref{lemuncb}, every well-defined multiplier $M_{m,\Phi,\Phi}$ with real symbol $m$ is self-adjoint.  Below we list some further  relations between the symbol and the operator type of a Bessel multiplier. In fact we investigate when a multiplier belongs to certain operator clones.

\begin{proposition} 
Let $M_{m,\Phi,\Psi}$ be a Bessel multiplier. 
 \begin{itemize} 
\item[{\rm (a)}] {\rm \cite{xxlmult1}} If $m\in c_0$, then  $M_{m,\Phi,\Psi}$ is a compact operator.
\item[{\rm (b)}] {\rm \cite{xxlmult1}} If $m\in \ell^1$, then  $M_{m,\Phi,\Psi}$ is a trace class operator with $\|M_{m,\Phi,\Psi}\|_{trace}\leq \sqrt{B_\Phi B_\Psi}\, \|m\|_1$ and ${\rm tr}(M_{m,\Phi,\Psi}) = \sum_n m_n \langle \phi_n,\psi_n\rangle$.
\item[{\rm (c)}] {\rm \cite{xxlmult1}} If $m\in \ell^2$, then  $M_{m,\Phi,\Psi}$ is a Hilbert Schmidt operator with $\|M_{m,\Phi,\Psi}\|_{HS}\leq \sqrt{B_\Phi B_\Psi}\,\|m\|_2$.
\item[{\rm (d)}]  If $m\in \ell^p$ for $1 < p < \infty$, then  $M_{m,\Phi,\Psi}$ is a Schatten-$p$ class operator with $\|M_{m,\Phi,\Psi}\|_{Sp}\leq \sqrt{B_\Phi B_\Psi}\, \|m\|_p$.
\end{itemize}
\end{proposition}
\bp
(d) follows directly from (b) and Prop. \ref{propgen1}(b) by complex interpolation \cite[Rem. 2.2.5, Theorems 2.2.6 and 2.2.7]{zhu1}.
\ep

\subsection{On invertibility} \label{subsinv}

As mentioned above, a multiplier is a time-variant filtering. As such it can be seen as the mathematical description of what a sound engineer does during a concert or recording session. Should the technician make an error, can we get the original signal back?
Or in the mathematical terms used here:
How and under which circumstances can we invert a multiplier?

To shorten notation, throughout this section $M$ denotes any one of the multipliers $M_{m,\Phi,\Psi}$ and $M_{m,\Psi,\Phi}$.

\subsubsection{Riesz multipliers, necessary and sufficient conditions for invertibility}

Schatten \cite{schatt1} investigated multipliers for orthonormal bases and showed many nice results leading to certain operator classes.
Extending the notion to Riesz bases keep the results very intuitive and easy, among them the following one:

\begin{proposition}\label{prr1}
Let $\Phi$ be a Riesz basis for $\mathcal{H}$. The following statements hold:
 \begin{itemize}
 \item[{\rm (a)}] {\rm \cite{xxlmult1}} 
  If $\Psi$ is a Riesz basis for $\mathcal{H}$ and  $m$ is semi-normalized, then 
 $M_{m,\Phi,\Psi}$ is invertible and 
 \begin{equation}\label{rminv}
M_{m,\Phi,\Psi}^{-1}=M_{1/m,\widetilde{\Psi},\widetilde{\Phi}}.
\end{equation}

\item[{\rm (b)}]  
{\rm \cite{balsto09new}}
 If $\Psi$ is a Riesz basis for $\mathcal{H}$, then $M$ is invertible if and only if  $m$ is semi-normalized.

\item[{\rm (c)}] {\rm \cite{balsto09new}} 
If $m$ is semi-normalized, then $M$ is invertible  if and only if $\Psi$ is a Riesz basis for $\mathcal{H}$. 

\end{itemize}
\end{proposition}

Under the assumptions of Proposition \ref{prr1}, one can easily observe that $M$ is invertible  if and only if $m\Psi$ is a Riesz basis for $\mathcal{H}$. The following proposition further clarifies the cases when $M$ is injective or surjective.

\begin{proposition} {\rm \cite{iwota11}}    Let $\Phi$ be a Riesz basis for $\mathcal{H}$. The following equivalences hold:

\begin{itemize}
\item[{\rm (a)}] 
 $M_{m,\Phi,\Psi}$ is  injective if and only if $m\Psi$ is a complete Bessel sequence in $\mathcal{H}$.

\item[{\rm (b)}] 
 $M_{m,\Psi,\Phi}$ is  injective if and only if $T_{m\Psi}$ is injective.

\item[{\rm (c)}] 
$M_{m,\Phi,\Psi}$ is  surjective  if and only if 
$\overline{m}\Psi$  is a 
Riesz sequence.
\item[{\rm (d)}] 
$M_{m,\Psi,\Phi}$ is surjective  if and only if $m\Psi$ is frame for $\mathcal{H}$.

\end{itemize}
\end{proposition}

For further results related to invertibility and non-invertibility of $M$ in the cases when $\Phi$ is a Riesz basis and $m$ is not necessarily semi-normalized, see \cite{iwota11}.

\subsubsection{Bessel multipliers, necessary conditions for invertibility}

Looking at invertible multipliers indicates somehow a kind of duality of the involved sequences. 
To make this more precise, let us state the following:

\begin{proposition}  {\rm \cite{balsto09new}}  Let $M_{m,\Phi,\Psi}$ be invertible.
\begin{itemize}
\item[{\rm (a)}] If $\Psi$ (resp. $\Phi$) is a Bessel sequence for $\mathcal{H}$ with bound $B$, then $m\Phi$ (resp. $m\Psi$)  satisfies the lower frame condition for $\mathcal{H}$ with bound $\frac{1}{B\, \|M_{m,\Phi,\Psi}^{-1}\|^2}$. 
\item[{\rm (b)}] If $\Psi$ (resp. $\Phi$) is a Bessel sequence in $\mathcal{H}$ and $m\in\ell^\infty$, then $\Phi$ (resp. $\Psi$) satisfies the lower frame condition for $\mathcal{H}$.
\item[{\rm (c)}]
If $\Psi$ and $\Phi$ are Bessel sequences in $\mathcal{H}$ and $m\in\ell^\infty$, then $\Psi$, $\Phi$, $m\Phi$, and $m\Psi$ are frames for $\mathcal{H}$.

\end{itemize}
\end{proposition}

\subsubsection{Frame multipliers}\label{invfrmult}

Naturally, the case for overcomplete frames is not that easy as for Riesz bases.  First we are interested to determine cases when multipliers for frames are invertible:

\vspace{.1in}
\noindent {\bf Sufficient conditions for invertibility and representation of the inverse via Neumann-like series}

\vspace{.1in} 

In this part of the section we present sufficient conditions for invertibility of multipliers $M_{m,\Phi,\Psi}$ based on perturbation conditions, and  formulas for the inverse $M_{m,\Phi,\Psi}^{-1}$ via Neumann-like series. 
In section \ref{sec-implement} we provide Matlab-codes for the inversion of multipliers according to Propositions \ref{mpos}, \ref{mpos2}, and \ref{p11}.

\vspace{.1in}
Let us begin our consideration with the specific case when $\Psi=\Phi$ and $m$ is positive (or negative) semi-normalized sequence. In this case the multiplier is simply a frame operator of an appropriate frame:
\begin{proposition} \label{psm} {\rm \cite[Lemma 4.4]{xxljpa1}} 
If $\Phi$ is a frame for $\mathcal{H}$ and $m$ is positive (resp. negative) and semi-normalized, then $M_{m,\Phi,\Phi}=S_{(\sqrt{m_n} \, \phi_n)}$  (resp. \sloppy\mbox{$M_{m,\Phi,\Phi}=-S_{(\sqrt{|m_n|} \, \phi_n)}$}) for the weighted frame $(\sqrt{m_n} \,\phi_n)$ and is therefore invertible on $\mathcal{H}$.
\end{proposition}

If we give up with the condition $\Phi=\Psi$, but still keep the condition on $m$ to be semi-normalized and positive (or negative), then the multiplier is not necessarily invertible and thus it is not necessarily representable as a frame operator. Consider for example the frames
$\Phi=(e_1, e_1, e_2, e_2, e_3, e_3, . . .)$ and $\Psi=(e1, e1, e2, e3, e4, e5, . . .)$; the multiplier $M_{(1),\Phi,\Psi}$ is
well-defined (even unconditionally convergent) but not injective \cite[Ex. 4.6.2]{bstable09}. 
For the class of multipliers which satisfy the assumptions of the above proposition except not necessarily the assumption $\Phi=\Psi$, below we present a sufficient condition for invertibility, which is a reformulation of \cite[Prop. 4.1]{balsto09new}. The reformulation is done at aiming efficiency of the implementation of the inversion of multiple multipliers - given $\Phi$ and $m$, once the frame operator of the weighted frame $(\sqrt{m_n}\phi_n)$ and its inverse are calculated, one can invert $M_{m,\Phi,\Psi}$ for different $\Psi$  fast,  just with operations of matrix summation and multiplication  
 (see Alg. \ref{alg:Prop8a} in Section \ref{sec-implement}).

 \begin{proposition} \label{mpos}

  Let $\Phi$ be a frame for $\mathcal{H}$, $m$ be a positive (or negative) semi-normalized sequence, and $a$ and $b$ satisfy $0< a \leq  \mid m_n \mid \leq b$ for every $n$. Assume that the sequence 
  $\Psi$ satisfies the condition 
 \begin{itemize}
 \item[$\ \ \ \ \ \mathcal{P}_1${\rm :} ]  
 \ \ \ \ \ \ \ \ \ \mbox{$\sum |\langle h, \psi_n-\phi_n \rangle|^2 \leq \mu \|h\|^2$, \ $\forall$ $h\in\mathcal{H}$,}
 \end{itemize}
 for some $\mu \in [0,  \frac{a^2A_\Phi^2}{b^2B_{\Phi}})$.
Then
 $\Psi$ is a frame for $\mathcal{H}$, 
$M$ is invertible on $\mathcal{H}$ and 
\begin{equation} \label{mnorm} 
\frac{1}{b B_\Phi + b\sqrt{\mu B_\Phi } } \|h\|\leq \|M^{-1} h\| \leq \frac{1}{aA_\Phi - b \sqrt{\mu B_\Phi } } \|h\|,
\end{equation}

$$M^{-1}=
\left\{\begin{array}{ll}
\sum_{k=0}^\infty [S_{(\sqrt{m_n}\phi_n)}^{-1}(S_{(\sqrt{m_n}\phi_n)}-M)]^k S_{(\sqrt{m_n}\phi_n)}^{-1}, \ &\mbox{if} \ m_n>0,\forall n,\\
\sum_{k=0}^\infty (-1)[S_{(\sqrt{|m_n|}\phi_n)}^{-1}(S_{(\sqrt{|m_n|}\phi_n)}+M)]^k S_{(\sqrt{|m_n|}\phi_n)}^{-1},  \ &\mbox{if} \ m_n<0, \forall n,
\end{array}
\right.$$

where the $n$-term error $\|M^{-1}-\sum_{k=0}^n \ldots\|$  
is bounded by 
\begin{equation}\label{errorbound}
\left( \frac{b \sqrt{\mu B_{\Phi} }}{a A_\Phi}\right) ^{n+1} \cdot
 \frac{1}{a A_\Phi-b \sqrt{\mu B_{\Phi}}}.
 \end{equation}
  \end{proposition}
 
\vspace{.1in} For $\mu=0$, the above statement gives Proposition \ref{psm}.
Note that $\mathcal{P}_1$ is a standard perturbation result for frames (for $\mu < A_\Phi$).
Notice that the bound $\frac{a^2A_\Phi^2}{b^2B_{\Phi}}$ for $\mu$ is sharp in the sense that when the inequality in $\mathcal{P}_1$ holds with $\mu\geq \frac{a^2A_\Phi^2}{b^2B_{\Phi}}$, the conclusions may fail (consider for example $m=(1)$, $\Phi=(e_n)_{n=1}^\infty$, $\Psi=(2e_1, \frac{1}{2}e_2, \frac{1}{3}e_3, \frac{1}{4}e_4,\ldots)$), but may also hold (consider for example $m=(1)$, $\Phi=(e_n)_{n=1}^\infty$, $\Psi=(2e_1, e_2, e_3, e_4,\ldots)$).

Note that for an overcomplete frame $\Phi$, a frame $\Psi$ satisfying $\mathcal{P}_1$ with $\mu < \frac{a^2A_\Phi^2}{b^2B_{\Phi}}(\leq A_\Phi)$  must also be overcomplete, by \cite[Cor. 22.1.5]{ole1}. This is also in correspondence with the aim to have an invertible frame multiplier in the above statement - when $m$ is semi-normalized and $\Phi$ and $\Psi$ are frames, then invertibility of $M_{m,\Phi,\Psi}$ requires either both $\Phi$ and $\Psi$ to be Riesz bases, or both to be overcomplete (see Prop. \ref{prr1}). 

Now we continue with consideration of more general cases, giving up the positivity/negativity of $m$, allowing  complex values. 
Again, aiming at an efficient inversion implementation for a set of multipliers (for given $\Phi$ and $m$, and varying $\Psi$), we reformulate the result from \cite[Prop. 4.1]{balsto09new}:

\begin{proposition} \label{mpos2}  
 Let $\Phi$ be a frame for $\mathcal{H}$ and 
let  $\lambda:=\sup_n |m_n-1| <\frac{A_\Phi}{B_\Phi}$. 
 Assume that the sequence $\Psi$ satisfies the condition 
 $\mathcal{P}_1$ with $\mu\in[0,\frac{(A_\Phi-\lambda B_\Phi)^2}{(\lambda+1)^2 B_\Phi})$.  
 Then $\Psi$ is a frame for $\mathcal{H}$, 
$M_{m,\Phi,\Phi}$ and $M$ are invertible on $\mathcal{H}$, and 
$$
\frac{1}{(\lambda+1) B_\Phi } \|h\|\leq \|M_{m,\Phi,\Phi}^{-1} h\| \leq 
\frac{1}{A_\Phi - \lambda B_\Phi  } \|h\|,
$$
$$
\frac{1}{(\lambda+1)( B_\Phi +\sqrt{\mu B_\Phi }) } \|h\|\leq \|M^{-1} h\| \leq 
\frac{1}{A_\Phi - \lambda B_\Phi -(\lambda+1) \sqrt{\mu B_\Phi } \,} \|h\|,
$$
\begin{equation}\label{mphiphi}
M_{m,\Phi,\Phi}^{-1}=
\sum_{k=0}^\infty [S_\Phi^{-1}(S_\Phi-M_{m,\Phi,\Phi})]^k S_\Phi^{-1}, 
\end{equation}
where the  $n$-term error is bounded by 
 $\left(\frac{\lambda B_{\Phi}}{A_\Phi}\right)^{n+1} \cdot
\frac{1}{A_\Phi - \lambda B_\Phi},$ and
\begin{equation}\label{mphipsi}
M^{-1}=
\sum_{k=0}^\infty [M_{m,\Phi,\Phi}^{-1}(M_{m,\Phi,\Phi}-M)]^k M_{m,\Phi,\Phi}^{-1}, 
\end{equation}
where the $n$-term error is bounded by 
 $\left(\frac{(\lambda+1)\sqrt{\mu B_{\Phi}}}{A_\Phi-\lambda B_\Phi}\right)^{n+1} \cdot
\frac{1}{A_\Phi - \lambda B_\Phi -(\lambda+1) \sqrt{\mu B_\Phi }}.$
 \end{proposition}

\vspace{.1in}
In the special case when $\Phi$ and $\Psi$ are dual frames, one has a  statement concerning invertibility of the multiplier $M_{m,\Phi,\Psi}$  with simpler formula (\ref{edno3}) for the inverse operator compare to the formula (\ref{mphipsi}):

\begin{proposition} {\rm \cite[Prop. 4.4]{balsto09new}} \label{p10}
Let $\Phi$ be a frame for $\mathcal{H}$ and let $\Psi$ be a dual frame of $\Phi$. 
Assume that $m$ is such that there is $\lambda$ satisfying $\left| m_n -1\right| \leq \lambda < \frac{1}{\sqrt{B_\Phi B_\Psi}}$ for all $n\in\mathbb{N}$. 
Then $M$ is invertible,  
 $$ \frac{1}{1+\lambda \sqrt{B_\Phi B_{\Phi^d}}}\|h\|\leq \|M^{-1} h\|\leq \frac{1}{1-\lambda \sqrt{B_\Phi B_{\Phi^d}}}\|h\|, \ \forall h\in \mathcal{H},$$ 
 \begin{equation}\label{edno3}  
 M^{-1}=\sum_{k=0}^\infty (I -M)^k,
 \end{equation}
 and the $n$-term error is bounded by 
$\frac{\left( \lambda \sqrt{B_\Phi B_\Psi} \right)^{n+1}}{1-\lambda \sqrt{B_\Phi B_\Psi}}$. 
\end{proposition}

Note that the bound  for $\lambda$ in the above proposition is sharp in the sense that one can not claim validity of the conclusions using a bigger constant instead of $\frac{1}{\sqrt{B_\Phi B_\Psi}}$. Indeed, if the assumptions hold with $\sup_n|m_n-1|=\lambda=\frac{1}{\sqrt{B_\Phi B_\Psi}}$, then the multiplier might not be invertible on $\mathcal{H}$, consider for example $M_{(1/n), (e_n), (e_n)}$ which is not surjective.

Here we extend Proposition \ref{p10} to approximate duals and provide implementation for the extended result. Recall that given a frame $\Phi$ for $\mathcal{H}$, the sequence $\Psi$ is called an \emph{approximate dual} 
 \emph{of  $\Phi$} \cite{chla08}  if for some $\varepsilon \in [0,1)$ one has $\|T_\Psi U_\Phi f-f\|\leq \varepsilon \|f\|$  and for all $f\in\mathcal{H}$ (in this case we will use the notion  \emph{$\varepsilon$-approximate dual of $\Phi$}).

\begin{proposition}  \label{p11}
Let $\Phi$ be a frame for $\mathcal{H}$ and let $\Psi$ be an $\varepsilon$-approximate dual frame of $\Phi$ ($\varepsilon\in [0,1)$).  
Assume that $m$ is such that $\lambda:=\sup_n\left| m_n -1\right|  < \frac{1-\varepsilon}{\sqrt{B_\Phi B_\Psi}}$ for all $n\in\mathbb{N}$. 
Then $M$ is invertible,  
$$  \frac{1}{1+\lambda \sqrt{B_\Phi B_{\Phi^d}}+\varepsilon}\|h\|\leq \|M^{-1} h\|\leq \frac{1}{1-\lambda \sqrt{B_\Phi B_{\Phi^d}}-\varepsilon}\|h\|, \ \forall h\in \mathcal{H},$$
 \begin{equation}\label{edno3a}  
 M^{-1}=\sum_{k=0}^\infty (I -M)^k,\end{equation}
 and the $n$-term error is bounded by 
$\frac{\left( \lambda \sqrt{B_\Phi B_\Psi} +\varepsilon\right)^{n+1}}{1-\lambda \sqrt{B_\Phi B_\Psi}-\varepsilon}$. 
\end{proposition}
\bp The case $\varepsilon=0$ gives Proposition \ref{p10}. For the more general case $\varepsilon\in [0,1)$, observe that for every $f\in\mathcal{H}$, 
$$
\|Mf-f\| \leq \|Mf - T_\Psi U_\Phi f\| + \| T_\Psi U_\Phi f -f\| 
\leq (\lambda \sqrt{B_\Phi B_\Psi} + \varepsilon) \|f\|.
$$ 
Applying \cite[Prop. 2.2]{balsto09new} (based on statements in \cite{GGKbook2003} and \cite{CCpert97}), one comes to the desired conclusions.
\ep

\vspace{.1in}
For the cases when $\Phi$ and $\Psi$ are equivalent frames, one
 can use the following sufficient conditions for invertibility and representations of the inverses:

\begin{proposition} {\rm \cite{balsto09new}} \label{gphi}
Let $\Phi$ be a frame for $\mathcal{H}$, $G:\mathcal{H}\to\mathcal{H}$ be a bounded bijective operator and $\psi_n=G\phi_n$, $\forall n$, i.e. $\Phi$ and $\Psi$ are equivalent frames. 
Let $m$ be semi-normalized and satisfy one of the following three conditions:
$m$ is positive; $m$ is negative; or
$\sup_n|m_n-1|< A_\Phi/B_\Phi$. 
Then $\Psi$ is a frame for $\mathcal{H}$, $M$ is invertible, 
$ M_{m,\Phi,\Psi}^{-1}=
(G^{-1})^* M^{-1}_{m,\Phi,\Phi}$ and 
 $M_{m,\Psi,\Phi}^{-1}=
 M^{-1}_{m,\Phi,\Phi} G^{-1}
$.
\end{proposition}

If one considers the canonical dual frame  $\widetilde{\Phi}$ of a frame $\Phi$, it is clearly equivalent to $\Phi$ and furthermore, it is the only dual frame of $\Phi$ which is equivalent to $\Phi$ \cite{hala00}. In this case one can apply both Propositions \ref{p10} and \ref{gphi}.

\vspace{.1in} \noindent 
{\bf Representation of the inverse as a multiplier using the reciprocal symbol and appropriate dual frames of the given ones}

\vspace{.1in}
Some of the results above were oriented to the representation of the inverse of an invertible frame multiplier via Neumann-like series. Here our attention is on  representation of the inverse  as a frame multiplier of a specific type. 
The inverse of any invertible frame multiplier $M_{m,\Phi,\Psi}$ with non-zero elements of $m$ can always be written using the reciprocal symbol $1/m$,  any given frame $G$ and an appropriate sequence related to a dual frame $(g_n^d)$ of $G$ (e.g. as $M_{m,\Phi,\Psi}^{-1}=M_{1/m, (M^{-1}(m_n g_n^d)), (g_n)}$) and actually, any invertible operator $V$ can be written as a multiplier in such a way 
($V=M_{1/m, (V(m_n g_n^d)), (g_n)}$), the focus here is on the possibility for representation of  $M_{m,\Phi,\Psi}^{-1}$ using the reciprocal symbol and  {\it specific dual frames} of $\Phi$ and $\Psi$, more precisely, using
$$M_{1/m, \Psi^d, \Phi^d}$$ for some dual frames $\Phi^d$ of $\Phi$  and $\Psi^d$  of $\Psi$.

The motivation for the consideration of such representations comes from Proposition \ref{prr1}(a), which concerns Riesz multipliers with semi-normalized weights and representation using the reciprocal symbol and the canonical duals (the only dual frames in this case) of the given Riesz bases. Since the representation  (\ref{rminv}) is not limited to Riesz multipliers  (as a simple example consider the case when $\Phi=\Psi$ is an overcomplete frame and $m=(1)$), and it is clearly not always valid, it naturally leads to the investigation of cases where such inversion formula holds for overcomplete frames. Furthermore, for the cases of non-validity of the formula, it opens the question whether the canonical duals can be replaced with other dual frames.

\begin{theorem} {\rm \cite{bsreprinv2015,SBSamptaExt16}} \label{ff3} 
Let $\Phi$ and $\Psi$ be frames for $\mathcal{H}$, and let 
the symbol  
$m$ be such that $m_n\neq 0$ for every $n$ and the sequence $m\Phi$ be a frame for $\mathcal{H}$. 
Assume that $M_{m, \Phi, \Psi}$ is invertible. 
Then the following statements hold.
\begin{itemize}
\item[{\rm (i)}]
There exists a unique sequence $\Psi^\dagger$ in $\mathcal{H}$ 
so that 
$$
 M_{m, \Phi, \Psi}^{-1}  = M_{1/m, \Psi^\dagger, \Phi^{ad}}, \ \forall \mbox{ a-pseudo-duals $\Phi^{ad}$  of $\Phi$,}
$$
 and it is a dual frame of $\Psi$. Furthermore, $\Psi^\dagger=(M_{m,\Phi,\Psi}^{-1}(m_n\phi_n))_{n=1}^\infty$.
 
 \item[{\rm (ii)}]
 If $G=(g_n)_{n=1}^\infty$ is a sequence (resp. Bessel sequence) in $\mathcal{H}$  such that $ M_{1/m, \Psi^\dagger, G}$ is well-defined and  
$$M_{m, \Phi, \Psi}^{-1}  = M_{1/m, \Psi^\dagger, G},$$
then $G$ must be an a-pseudo-dual (resp. dual frame) of $\Phi$.

\end{itemize}
\end{theorem} 

Note that the uniqueness of $\Psi^\dagger$ in the above theorem is even guaranteed from the validity of $ M_{m, \Phi, \Psi}^{-1}  = M_{1/m, \Psi^\dagger, \Phi^{d}}$ for all dual frames $\Phi^{d}$  of $\Phi$. Indeed, take $\Psi^\dagger$ to be determined by the above theorem and assume that there is a sequence $F$ such that
$M_{m, \Phi, \Psi}^{-1}  = M_{1/m, F, \Phi^{d}}$ for all dual frames $\Phi^{d}$  of $\Phi$. Then $M_{1/m, F-\Psi^\dagger, \Phi^{d}}$ for all dual frames $\Phi^{d}$  of $\Phi$, which by \cite[Lemma 3.2]{SBSamptaExt16} implies that $F=\Psi^\dagger$. 

\vspace{.1in}
Similar statements hold with respect to appropriate dual frame of $\Phi$:

\begin{theorem} \label{ff32} {\rm \cite{bsreprinv2015,SBSamptaExt16}} 
Let $\Phi$ and $\Psi$ be frames for $\mathcal{H}$, and let 
the symbol  
$m$ be such that $m_n\neq 0$ for every $n$ and the sequence $m\Psi$ is a frame for $\mathcal{H}$. 
Assume that $M_{m, \Phi, \Psi}$ is invertible. 
Then the following statements hold.
\begin{itemize}
\item[{\rm (i)}]
There exists a unique sequence $\Phi^\dagger$ in $\mathcal{H}$ 
so that 
$$
 M_{m, \Phi, \Psi}^{-1}  = M_{1/m, \Psi^{sd}, \Phi^\dagger}, \ \forall \mbox{ s-pseudo-duals $\Psi^{sd}$ of $\Psi$,}
$$
 and it is a dual frame of $\Phi$.  Furthermore, $\Phi^\dagger=((M_{m,\Phi,\Psi}^{-1})^*(\overline{m_n}\psi_n))_{n=1}^\infty$.
\item[{\rm (ii)}]
If $F=(f_n)_{n=1}^\infty$ is a sequence (resp. Bessel sequence) in $\mathcal{H}$   such that $M_{1/m, F, \Phi^\dagger}$is well-defined and
 $$M_{m, \Phi, \Psi}^{-1}  = M_{1/m, F, \Phi^\dagger},$$
then $F$ must be an s-pseudo-dual (resp. dual frame)  of  $\Psi$.
\end{itemize}
\end{theorem} 

   \begin{corollary} \label{ff4} {\rm \cite{SBSamptaExt16}} 
Let $\Phi$ and $\Psi$ be frames for $\mathcal{H}$, and let 
the symbol  
$m$ be such that $m_n\neq 0$ for every $n$ and $m\in\ell^\infty$.  
Assume that $M_{m, \Phi, \Psi}$ is invertible. 
Then Theorems \ref{ff3} and  \ref{ff32} apply. 
\end{corollary}

Below we consider the question for representation of the inverse using the canonical dual frames of the given ones and the question when the special dual $\Phi^\dagger$ (resp. $\Psi^\dagger$) coincides with the canonical dual of $\Phi$ (resp. $\Psi$).

\begin{proposition} \label{q} {\rm \cite{bsreprinv2015,SBSamptaExt16}}
Let $\Phi$ and $\Psi$ be frames for ${\mathcal{H}}$ and let $m$ be semi-normalized. Assume that $M_{m,\Phi,\Psi}$ be invertible.   
Then 
\begin{equation} \label{psid}
 M_{m,\Phi,\Psi}^{-1}=M_{1/m,\widetilde{\Psi},\widetilde{\Phi}} \Leftarrow 
 \Psi \mbox{ is equivalent   to } m\Phi 
\Leftrightarrow
 \Psi^\dagger = \widetilde{\Psi};
\end{equation}
\begin{equation} \label{phid}
M_{m,\Phi,\Psi}^{-1}=M_{1/m,\widetilde{\Psi},\widetilde{\Phi}} \Leftarrow 
\Phi \mbox{ is equivalent   to } \overline{m}\Psi
\Leftrightarrow 
\Phi^\dagger = \widetilde{\Phi}.
\end{equation}
$$M_{m,\Phi,\Psi}^{-1}=M_{1/m,\widetilde{\Psi},\widetilde{\Phi}} \nRightarrow 
(\Phi \mbox{ is equivalent   to } \overline{m}\Psi)
\mbox{ or }
  (\Psi \mbox{ is equivalent   to } m\Phi). $$
\vspace{.07in}
\noindent 
\end{proposition}

In the cases when the symbol $m$ is a constant sequence, the one-way implications in the above proposition become equivalences, more precisely, the following holds:

\begin{proposition} \label{qconstant} {\rm \cite{bsreprinv2015}}
Let $\Phi$ and $\Psi$ be frames for ${\mathcal{H}}$. 
For  symbols of the form $m=(c,c,c,\ldots)$, $c\neq 0$,  the following equivalences hold:
\begin{eqnarray*}
M_{m,\Phi,\Psi} \mbox{ is invertible and }  M_{(c),\Phi,\Psi}^{-1}=M_{(1/c),\widetilde{\Psi},\widetilde{\Phi}} 
&\Leftrightarrow& 
\Psi \mbox{ is equivalent   to } \Phi \\
&\Leftrightarrow&  M_{m,\Phi,\Psi} \mbox{ is invertible and }  \Psi^\dagger = \widetilde{\Psi} 
\\
&\Leftrightarrow& M_{m,\Phi,\Psi} \mbox{ is invertible and }  \Phi^\dagger = \widetilde{\Phi}.
\end{eqnarray*} 
\end{proposition}

Motivated by Theorems \ref{ff3} and \ref{ff32}, it is natural to consider the question whether (\ref{psid}) and (\ref{phid}) from Proposition \ref{q} hold under more general assumptions. The answer is "positive":

\begin{proposition} \label{p1}
Let $\Phi$ and $\Psi$ be frames for $\mathcal{H}$, and let 
the symbol  
$m$ be such that $m_n\neq 0$ for every $n$ and the sequence $m\Phi$ be a frame for $\mathcal{H}$. 
Then the following holds. 
\begin{eqnarray*}
M_{m,\Phi,\Psi} \mbox{ is invertible and }  \Psi^\dagger = \widetilde{\Psi} 
& \Leftrightarrow&  \Psi \mbox{ is equivalent   to } m\Phi \\
& \Rightarrow& M_{m,\Phi,\Psi} \mbox{ is invertible and } M_{m,\Phi,\Psi}^{-1}=M_{1/m,\widetilde{\Psi},\widetilde{\Phi}}
\end{eqnarray*}
\end{proposition}
\bp
First let $M_{m,\Phi,\Psi}$ be invertible and $\Psi^\dagger=\widetilde{\Psi}$, where $\Psi^\dagger$ is determined based on Theorem \ref{ff3}. Then using Theorem \ref{ff3}(i), one obtains $M_{m,\Phi,\Psi}^{-1}=M_{1/m,\widetilde{\Psi},\widetilde{\Phi}}$ and 
furthermore 
  $\psi_n=S_{\widetilde{\Psi}}^{-1} (\widetilde{\psi_n})=
  S_{\widetilde{\Psi}}^{-1} (\psi_n^\dagger)=
  S_{\widetilde{\Psi}}^{-1} M_{m,\Phi,\Psi}^{-1}(m_n\phi_n)$, $\forall n$, leading to equivalence of the frames $\Psi$ and $m\Phi$.

Conversely, let the frames $\Psi$ and $m\Phi$ be equivalent. Then by Proposition \ref{qconstant}, the multiplier $M_{(1), m\Phi,\Psi}(=M_{m,\Phi,\Psi}$) is invertible, and thus, by Theorem \ref{ff3}, a dual frame $\Psi^\dagger$ of $\Psi$ is determined by  $\Psi^\dagger=(M_{m,\Phi,\Psi}^{-1}(m_n\phi_n))_{n=1}^\infty$.
This implies that the dual frame $\Psi^\dagger$ of $\Psi$ is at the same time equivalent to the frame $\Psi$, which by \cite[Prop. 1.14]{hala00}
 implies that  $\Psi^\dagger$ must be the canonical dual of $\Psi$. 
\ep

 In a similar way as above, one can state a corresponding result involving $\Phi^\dagger$:
\begin{proposition} \label{p2}
Let $\Phi$ and $\Psi$ be frames for $\mathcal{H}$, and let 
the symbol  
$m$ be such that $m_n\neq 0$ for every $n$ and the sequence $m\Psi$ is a frame for $\mathcal{H}$. Then the following holds. 
\begin{eqnarray*} M_{m,\Phi,\Psi} \mbox{ is invertible and }  \Phi^\dagger = \widetilde{\Phi} 
& \Leftrightarrow&  \Phi \mbox{ is equivalent   to } \overline{m}\Psi \\
& \Rightarrow& M_{m,\Phi,\Psi} \mbox{ is invertible and } M_{m,\Phi,\Psi}^{-1}=M_{1/m,\widetilde{\Psi},\widetilde{\Phi}}
\end{eqnarray*}

\end{proposition}

For examples, which illustrate statements from this section, see \cite{bsreprinv2015,SBSamptaExt16}.

\section{Time-Frequency Multipliers}\label{Sec:exampl}

In this section we focus on sequences with particular structure, namely, on Gabor and wavelet sequences, which are very important for applications.  
Considering the topic of unconditional convergence, we apply results from Section \ref{subs-UncondConv} directly, while for the topic of invertibility we go beyond the presented results in Section \ref{subsinv} and consider further questions motivated by the specific structure of the sequences.

\subsection*{On the unconditional convergence}

As a consequence of Proposition \ref{propgen1}, for Gabor and wavelet systems, 
we have the following statements:

\begin{corollary}
Let $M_{m,\Phi,\Psi}$ be a Gabor (resp. wavelet) multiplier.
\begin{itemize}
\item[{\rm (c)}]
If $M_{m,\Phi,\Psi}$ is furthermore a Bessel multiplier, then it is unconditionally convergent if and only if $m\in\ell^\infty$.
\item[{\rm (d)}]
If $m$ is NBB, then $M_{m,\Phi,\Psi}$ is unconditionally convergent if and only if  
 $M_{m,\Phi,\Psi}$ is a Bessel multiplier and $m$ is semi-normalized.
  \item[{\rm (e)}]
  If $\Phi$ is Bessel (resp. $\Phi$ is Bessel and $m$ is NBB), then  $M_{m,\Phi,\Psi}$ is unconditionally convergent if and only if
$m\Psi$ is Bessel (resp. $\Psi$ is Bessel).
\item[{\rm (f)}] 
If $M_{m,\Phi,\Psi}$  is furthermore a Riesz multiplier, then it is well-defined  if and only if it is unconditionally convergent if and only if $m\in\ell^\infty$. 
\end{itemize}

\end{corollary}

\subsection*{On the invertibility}

For the particular case of Gabor and wavelet multipliers, one can apply the general  invertibility results from Section \ref{subsinv} and in addition, one can be interested in cases where the dual frames $\Phi^\dagger$, $\Psi^\dagger$, induced by an invertible Gabor (resp. wavelet) frame multiplier have also a Gabor (resp. wavelet) structure and whether one can write the inverse multiplier as a Gabor (resp. wavelet) multiplier.

\begin{proposition}
 \label{corgabor2} {\rm \cite{sbsampta17}} 
Let $\Lambda=\{ (\omega, \tau)\}$  be a lattice in ${\mathbb R}^{2d}$, i.e., a discrete subgroup of ${\mathbb R}^{2d}$ of the form $A{\mathbb Z}^{2d}$ for some invertible matrix $A$, and let $\Phi=(E_\omega T_\tau v)_{(\omega,\tau)\in\Lambda}$ and $\Psi=(E_\omega T_\tau u)_{(\omega,\tau)\in\Lambda}$  be (Gabor) frames. 
Assume that the (Gabor) frame type operator $V=M_{(1),\Phi,\Psi}$ is invertible. Then 
the following holds.
\begin{itemize}
\item[{\rm (a)}] The frames $\Phi^\dagger$ and $\Psi^\dagger$ determined by Theorems \ref{ff3} and \ref{ff32} are Gabor frames.
\item[{\rm (b)}]  $V^{-1}$ can be written as a  Gabor frame-type operator as follows: \\ 
$$ V^{-1} = M_{(1),(E_\omega T_\tau V^{-1}v )_{(\omega,\tau)\in\Lambda},\Phi^d} 
=  M_{(1),\Psi^d,(E_\omega T_\tau (V^{-1})^*u)_{(\omega,\tau)\in\Lambda}},
$$ 
using any dual Gabor frames  $\Phi^d$ and $\Psi^d$ of $\Phi$ and $\Psi$, respectively (in particular, the canonical duals). 

\item[{\rm (c)}] Let $(E_\omega T_\tau g)_{(\omega,\tau)\in\Lambda}$ be a Gabor frame for $L^2({\mathbb R})$. Then $V^{-1}$
can be written as the Gabor frame type operator 
$$M_{(1),(E_\omega T_\tau g)_{(\omega,\tau)\in\Lambda},({\widetilde h}_{\omega,\tau})_{(\omega,\tau)\in\Lambda}},$$ where $h_{\omega,\tau}=E_\omega T_\tau Vg$, $(\omega,\tau)\in\Lambda$.
\end{itemize}
\end{proposition}

 As an illustration of the above proposition, one may consider the following example with Gabor frame-type operators:
\begin{example} {\rm \cite{sbsampta17}}
Let $\Psi=(E_{kb}T_{na} u)_{k,n\in\mathbb{Z}}$ be a frame for $L^2(\mathbb{R})$ (for example, take $u(x)$ to be the Gaussian $e^{-x^2}$ and $a\in\mathbb{R}$ and $b\in\mathbb{R}$ so that $ab<1$). Furthermore, let $F$ be an invertible operator on $L^2(\mathbb{R})$ which commutes with all  $E_{kb}T_{na}$, $k,n\in\mathbb{Z}$, (e.g., take $F=E_{k_0/a}T_{n_0/b}$ for some $k_0,n_0\in\mathbb{Z}$) and let 
$\Phi=(E_{kb}T_{na} Fu)_{k,n\in\mathbb{Z}}$. 
Then the Gabor frame-type operator $V=M_{(1),\Phi,\Psi}(=FS_\Psi)$ is invertible and furthermore:
\begin{itemize}
\item[{\rm (a)}] $\Phi^\dagger=(E_{kb}T_{na} (V^{-1})^*u)_{k,n\in\mathbb{Z}}=\widetilde{\Phi}$ and
$\Psi^\dagger=(E_{kb}T_{na} S_\Psi^{-1}u)_{k,n\in\mathbb{Z}}=\widetilde{\Psi}$; 
\item[{\rm (b)}] 
$V^{-1} = M_{(1),\psi^\dagger,\Phi^d} = M_{(1),\widetilde{\Psi},\Phi^d}$
for any dual Gabor frame  $\Phi^d$ of   $\Phi$. In particular, chosing the canonical dual $\widetilde{\Phi}$, one comes to the formula $V^{-1} = 
M_{(1),\widetilde{\Psi},\widetilde{\Phi}}$, whose validity is in correspondence with Proposition \ref{qconstant};

\item[{\rm (c)}] 
applying Proposition \ref{corgabor2}(c) with $g=u$, we have 
$V^{-1} = M_{(1),\Psi,(\widetilde{h}_{k,n})_{k,n\in\mathbb{Z}}}$, where $h_{k,n}=
E_{kb}T_{na}Vu$;

applying Proposition \ref{corgabor2}(c) with $g=S^{-1}u$, we obtain
$V^{-1} = M_{(1),\widetilde{\Psi},\widetilde{\Phi}}$.
\end{itemize}
\end{example}

The representation of the inverse of a Gabor frame type operator as a Gabor frame type operator has turned out to be related to commutative properties of the mutliplier with the time-frequency shifts:

\begin{proposition} \label{general1} {\rm \cite{bsreprinv2015}}
Let $\Lambda=\{ (\omega, \tau)\}$  be a lattice in ${\mathbb R}^{2d}$, i.e., a discrete subgroup of ${\mathbb R}^{2d}$ of the form $A{\mathbb Z}^{2d}$ for some invertible matrix $A$, 
$g\in L^2(\mathbb{R}^d)$, and  $(E_\omega T_\tau  g)_{(\omega,\tau)\in\Lambda}$ be a Gabor frame for $ L^2(\mathbb{R}^d)$. Let $V:  L^2(\mathbb{R}^d)\to  L^2(\mathbb{R}^d)$ be a bounded bijective operator. Then 
the following statements are equivalent.

\begin{itemize}
\item[{\rm (a)}] For every $(\omega, \tau)\in\Lambda$, $V E_\omega T_\tau g =E_\omega T_\tau Vg$.
\item[{\rm (b)}] For every $(\omega, \tau)\in\Lambda$ and every  $f\in L^2(\mathbb{R}^d)$, $V E_\omega T_\tau f = E_\omega T_\tau V f$ (i.e.,  $V$ commutes with $E_\omega T_\tau$ for every $(\omega, \tau)\in\Lambda$).

\item[{\rm (c)}]  $V$ can be written as a Gabor frame-type multiplier with the constant symbol $(1)$ and with respect to the lattice $\Lambda$, i.e., $V$ is a Gabor frame-type operator with respect to the lattice $\Lambda$.

\item[{\rm (d)}]
$V^{-1}$ can be written as a Gabor frame-type multiplier with the constant symbol $(1)$ and with respect to the lattice $\Lambda$, i.e., $V$ is a Gabor frame-type operator with respect to the lattice $\Lambda$.
\end{itemize}

\end{proposition}

The above two propositions \ref{corgabor2} and \ref{general1}  concern the cases when the symbol of an invertible Gabor frame multiplier is the constant sequence $(1)$. They show that such multipliers commute with the corresponding time-frequencsy shifts and the frames induced by these multipliers have Gabor structure. 
 Considering the more general case of $m$ not necessarily being a constant sequence,  the commutative property of the multiplier with the time-frequency shifts $E_\omega T_\tau$ is no longer necessarily valid (see \cite[Corol. 5.2 and Ex. 5.3 ]{sbsampta17}), but it can still happens that the frames $(E_\omega T_\tau v)_{(\omega,\tau)\in\Lambda}^\dagger$ and $(E_\omega T_\tau u)_{(\omega,\tau)\in\Lambda}^\dagger$ determined by Theorems \ref{ff3} and \ref{ff32} have Gabor structure (see \cite[Ex. 5.6]{sbsampta17}). 
Below we consider cases, when validity of the commutativity and the desire to have $\Phi^\dagger$ (resp. $\Psi^\dagger$) being Gabor frames can only happen when $m$ is a constant sequence:

\begin{proposition}\label{corconstant} Let $\Lambda=\{ (\omega, \tau)\}$  be a lattice in ${\mathbb R}^{2d}$, $\Phi =(E_{\omega} T_{\tau} v)_{(\omega,\tau)\in\Lambda}$ 
and $\Psi=(E_{\omega} T_{\tau} u)_{(\omega,\tau)\in\Lambda}$ 
be Gabor frames for $L^2({\mathbb R})$, and 
$m=(m_{\omega,\tau})_{(\omega,\tau)\in\Lambda}$ be such that $m_{\omega,\tau}\neq 0$ for every $(\omega,\tau)\in\Lambda$.  
Assume that the Gabor frame multiplier $M_{m,\Phi,\Psi}$ is invertible on $L^2({\mathbb R})$ 
and  commutes with $E_{\omega} T_{\tau}$  for every $(\omega,\tau)\in\Lambda$. 
Then the following statements hold.
\begin{itemize}
\item[{\rm (a)}] \cite[Coroll. 5.5]{sbsampta17} If $m\Phi$ is a frame for $L^2({\mathbb R})$ and the frame $\Psi^\dagger$  determined by Theorem  \ref{ff3} has  Gabor structure, then $m$ must be a constant sequence. 

 \item[{\rm (b)}] \cite[Coroll. 5.5]{sbsampta17} If  $m\Psi$ is a frame for $L^2({\mathbb R})$ and the frame $\Phi^\dagger$  determined by Theorem  \ref{ff3} has  Gabor structure, then $m$ must be a constant sequence. 
 \item[{\rm (c)}] If $\Phi$ is a Riesz sequence, then $m$ must be a constant sequence. 
\end{itemize}
\end{proposition}
\bp
(c) 
 Fix an arbitrary couple $(\omega',\tau')\in\Lambda$. Having in mind the general relation 
\begin{eqnarray*}\label{generalrel}
& & M_{(m_{\omega,\tau})_{(\omega,\tau)\in\Lambda},\Phi,\Psi} E_{\omega'} T_{\tau'} =
E_{\omega'} T_{\tau'}   M_{(m_{\omega+\omega', \tau+\tau'})_{(\omega,\tau)\in\Lambda},\Phi,\Psi} 
\end{eqnarray*} 
(see \cite[Coroll. 5.2]{sbsampta17}), and using the assumption that $M_{m,\Phi,\Psi}$ commutes with $E_{\omega'} T_{\tau'} $ and the invertibility of $E_{\omega'} T_{\tau'} $, 
 we get 
$$M_{(m_{\omega,\tau})_{(\omega,\tau)\in\Lambda}, \Phi,\Psi}
=  M_{(m_{\omega+\omega', \tau+\tau'})_{(\omega,\tau)\in\Lambda},\Phi,\Psi}.
$$ 
The assumption that $\Phi$ is a Riesz sequence leads now to the conclusion that
$$ (m_{\omega, \tau}-m_{\omega+\omega', \tau+\tau'})  \langle f, E_\omega T_\tau u\rangle = 0, \ \forall (\omega,\tau)\in\Lambda, \ \forall f\in \mathcal{H}.$$
Since $u\neq 0$ (as otherwise $M_{m,\Phi,\Psi}$ is not invertible), one can apply $f=E_\omega T_\tau u\neq 0$ above and conclude that 
$m_{\omega, \tau}=m_{\omega+\omega', \tau+\tau'}$ for every $(\omega,\tau)\in\Lambda.$ 
Then, for any fixed  $(\omega,\tau)\in\Lambda$, we have 
$$m_{\omega, \tau}=m_{\omega+\omega', \tau+\tau'}, \  \forall (\omega',\tau')\in\Lambda,$$
implying that $m$ is a constant sequence. 
\ep

Observe that using the language of tensor products,   in a similar way as in  Proposition \ref{corconstant}(c)  one can show that  $\phi_\lambda \otimes \psi_\lambda$ being a Riesz sequence in the tensor product space $ \mathcal{H} \otimes \mathcal{H}$ would lead to  the conclusion that $m$ must be a constant sequence.

\vspace{.1in}
Notice that when m is not a constant sequence, it is still possible to have $\Phi^\dagger$ and $\Psi^\dagger$ with the Gabor structure (see \cite[Ex. 5.6]{sbsampta17}), which opens interest to further investigation of this topic.

In this manuscript we presented results focusing on Gabor frames.
Note that some of the above results hold for a more general class of frames, containing Gabor frames and coherent frames (so called pseudo-coherent frames), {defined  in \cite{bsreprinv2015} and further investigated in \cite{sbsampta17}. }
Note that the class of wavelet frames is not included in this class and in general, one can not state corresponding results like Propositions  \ref{corgabor2} and \ref{general1} for wavelet frames. The wavelet case is a topic of interest for further investigations of questions related to the structure of $\Phi^\dagger$ and $\Psi^\dagger$ for wavelet multipliers. This might be related to the topic of localized frames, see e.g. \cite{xxlgro14}, 
where wavelets also cannot be included.

\section{Implementation of the inversion of multipliers according to Section \ref{invfrmult}} \label{sec-implement}

For the inversion of 
multipliers $M_{m,\Phi,\Psi}$ and $M_{m,\Psi,\Phi}$ according to  Propositions \ref{mpos}, \ref{mpos2}, and \ref{p11}, we provide Matlab-LTFAT codes (including demo-files), which are available on the  webpage
\url{https://www.kfs.oeaw.ac.at/InversionOfFrameMultipliers}. The inversion is done
 using an iterative process based on the formulas from the corresponding statement.
Here we present the algorithms for these codes. Note that these implementations are only meant as proof-of-concept and could be made significantly more efficient, which will be the topic of future work, but is beyond the scope of this manuscript.

\begin{figure}[!h] 
\begin{center}
\includegraphics[width=0.65\textwidth]{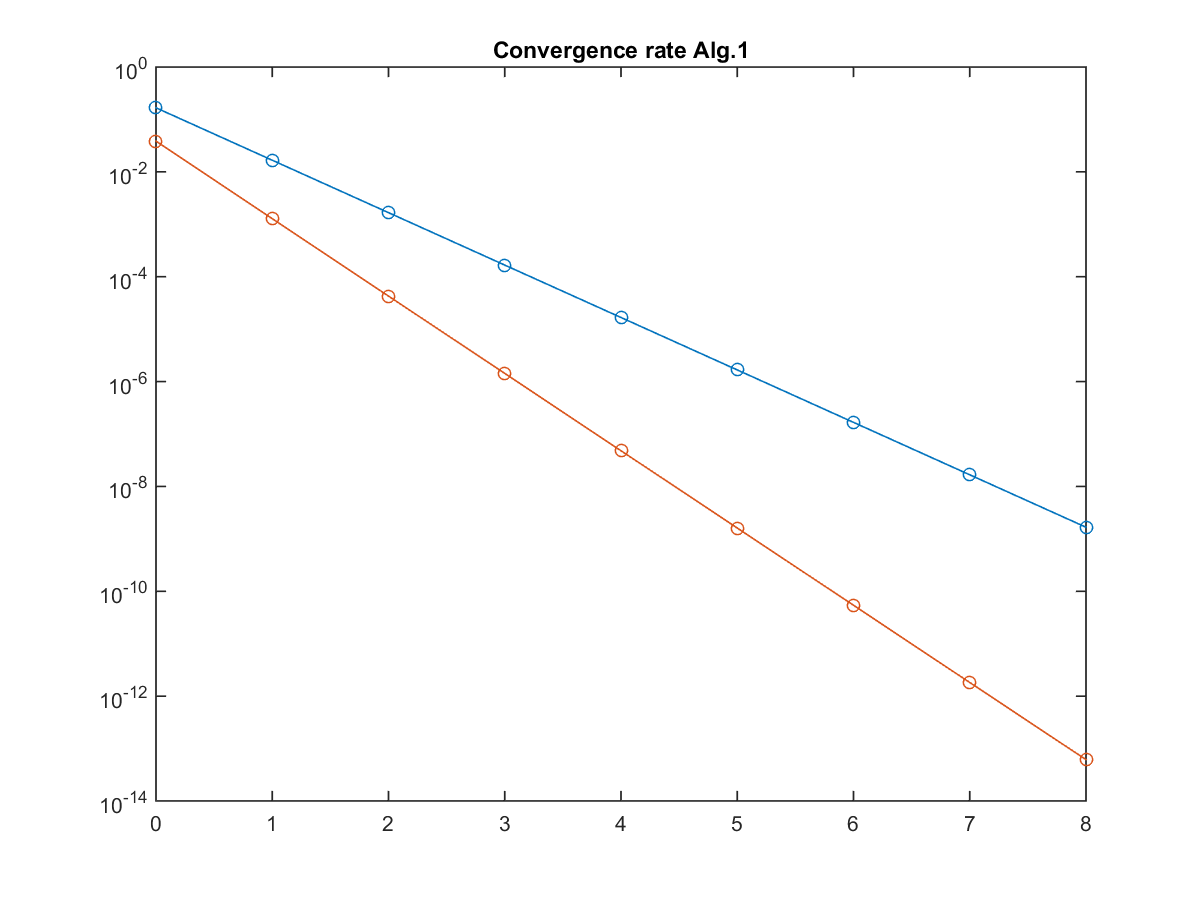} 
\caption{\label{fig:ConvRateGab} \small {\em  The convergence rate of Alg. \ref{alg:Prop8c} using base-$10$ logarithmic scale in the vertical axis and a linear scale in the horizontal axis. Here the absolute error in each iteration is plotted in red, and the convergence value predicted in Proposition \ref{mpos} is plotted in blue.
}} 
\end{center}
\end{figure}

As an implementation of Proposition \ref{mpos}, we provide three Matlab-LTFAT codes  - one for 
 computing $M_{m,\Phi,\Psi}^{-1}$ (see Alg. \ref{alg:Prop8a}), a second one for computing $M_{m,\Phi,\Psi}^{-1}f$ for given $f$ (see Alg. \ref{alg:Prop8b}), and a third one - for computing $M_{m,\Phi,\Psi}^{-1}$ in the case when $\Phi$ and $\Psi$ are Gabor frames (see Alg. \ref{alg:Prop8c}). For illustration of the convergence rate of  Alg. \ref{alg:Prop8c}, see Figure \ref{fig:ConvRateGab}. On the horizontal axis one has the number of the iterations in a linear scale, and on the vertical axes - the absolute error (the norm of the difference between the $n$-th iteration and the real inverse) in a logarithmic scale.
To produce this plot we have used Gabor frames with length of the transform $L= 4096$, time-shift $a= 1024$, number of channels $M= 2048$, a Hann window for $\Phi$, a Gaussian window for $\Psi$, symbol $m$ with elements between $1/2$ and $1$, and parameter $e=10^{-8}$, which guides the number of iterations by the predicted error bound (\ref{errorbound}). 
The size of the multiplier-matrix in this test is $4096 \times 4096$. 
The inversion takes only a few iterations - it stops at $n=8$ with error 
$6.1514 \cdot 10^{-14}$.

We also provide implementations of Prop.  \ref{mpos2} (see Alg. \ref{alg:Prop9}) and Prop.  \ref{p11} (see Alg. \ref{alg:Prop11}).

Note that the codes provide the possibility for multiple inversion of multipliers varying the frame $\Psi$. Although the initial step in the implementations of Propositions \ref{mpos} and \ref{mpos2} involves inversion of an appropriate frame operator and is computationally demanding, it depends only on $\Phi$ and $m$, and thus multiple inversions with different $\Psi$ will, in total, be very fast. 
It is the aim of our next work to improve the above mentioned algorithms avoiding the inversion in the initial step.

Notice that the iterative inversion according to  Prop.  \ref{p11} seems to be very promising for fast inversions as it does not involve an inversion of a full matrix.

In further work we will run tests to compare the speed of the inversion of the above implementation algorithms with inversion algorithms in Matlab and LTFAT, with the goal to improve the numerical efficiency.

\begin{algorithm}  \small
	\caption{Iterative computation of $M_{m,\Phi,\Psi}^{-1}$ (M1inv):  [TPsi,M1,M2,M1inv,M2inv,n] = \textbf{Prop8InvMultOp}(c,r,TPhi,TG,m,e)	}
	\label{alg:Prop8a}
	\begin{algorithmic}[1]
	    \State $T_G\gets T_G, m, T_\Phi$
	    \State $T_\Psi=T_\Phi+T_G, M1= T_\Phi*diag(m)*T_\Psi'$
		\State Initialize  $M1inv=S^{-1}_{(\sqrt{m_n}\phi_n)}$
		\State $n\gets e$, $m$, $T_\Phi$, $T_G$		
		\State {\bf if $n>0$}
		\State \ \ \  $P1=S^{-1}_{(\sqrt{m_n}\phi_n)}*(S_{(\sqrt{m_n}\phi_n)}-M1)$
		\State  \ \ \  Initialize $Q1=S^{-1}_{(\sqrt{m_n}\phi_n)}$
		\State  \ \ \ {\bf for i:=1 to n}
		\State  \ \ \ \ \ \ \ $Q1=P1*Q1$
		\State  \ \ \ \ \ \ \  $M1inv=M1inv + Q1$
		\State  \ \ \ {\bf end}
		\State {\bf end}
		\State $M_{m,\Phi,\Psi}^{-1} \gets M1inv$
	\end{algorithmic}
\end{algorithm}

\begin{algorithm} \small
	\caption{Iterative computation of $M_{m,\Phi,\Psi}^{-1}f$ (M1invf):  [TPsi,M1,M2,M1invf,M2invf,n] = \textbf{Prop8InvMultf}(c,r,TPhi,TG,m,f,e)}
	\label{alg:Prop8b}
	\begin{algorithmic}[1]
	    \State $T_G\gets T_G$, $m$, $T_\Phi$
	    \State $T_\Psi=T_\Phi+T_G$, $M1= T_\Phi*diag(m)*T_\Psi'$
		\State Initialize  $M1invf=S^{-1}_{(\sqrt{m_n}\phi_n)}f$
		\State $n\gets e$,  $m$, $T_\Phi$, $T_G$		
		\State {\bf if $n>0$}
		\State  \ \ \  $P1=S^{-1}_{(\sqrt{m_n}\phi_n)}*(S_{(\sqrt{m_n}\phi_n)}-M1)$
		\State  \ \ \  Initialize $q1=M1invf$
		\State \ \ \  {\bf for i:=1 to n}
		\State \ \ \  \ \ \  \ $q1=P1*q1$
		\State \ \ \  \ \ \  \ $M1invf=M1invf + q1$
		\State \ \ \  {\bf end}
		\State {\bf end}
		\State $M_{m,\Phi,\Psi}^{-1}f \gets M1invf$
	\end{algorithmic}
\end{algorithm}

\begin{algorithm}  \small
	\caption{Iterative computation of $M_{m,\Phi,\Psi}^{-1}$ (M1inv) for Gabor frames:  [TPhi,TPsi,M1,M2,M1inv,M2inv,n] = \textbf{Prop8InvMultOpGabor}(L,a,M,gPhi,gG,m,e)	}
	\label{alg:Prop8c}
	\begin{algorithmic}[1]
	   \State $\Phi\gets gPhi,a,M$
		   \State $T_\Phi\gets \Phi, L$

	  \State $G\gets gG,a,M$
		   \State $T_G\gets G, L$	
		 \State-17: These steps are like 1-13 of Alg. \ref{alg:Prop8a}.
		
	\end{algorithmic}
\end{algorithm}

\begin{algorithm} \small
	\caption{Iterative computation of $M_{m,\Phi,\Phi}^{-1}$ (M0inv) and $M_{m,\Phi,\Psi}^{-1}$ (M1inv):  [m,TPsi,M0,M1,M2,M0inv,n0,M1inv,M2inv,n] = \textbf{Prop9InvMultOp}(c,r,TPhi,TG,m,e)	}
	\label{alg:Prop9}
	\begin{algorithmic}[1]
	    \State $m\gets m$, $T_\Phi$
	    \State  $M0= T_\Phi*diag(m)*T_\Phi'$
	    \State Initialize  $M0inv=S_\Phi^{-1}$
		\State $n0\gets e$, $m$, $T_\Phi$		
		\State {\bf if $n0>0$}
		\State  \ \ \  $P0=S_\Phi^{-1}*(S_\Phi-M0)$
		\State \ \ \  Initialize $Q0=S_\Phi^{-1}$
		\State \ \ \  {\bf for i:=1 to n0}
		\State \ \ \  \ \ \  \ $Q0=P0*Q0$
		\State \ \ \  \ \ \  \ $M0inv=M0inv + Q0$
		\State \ \ \  {\bf end}
		\State {\bf end}
		\State $M_{m,\Phi,\Phi}^{-1} \gets M0inv	$		
		\State $T_G\gets T_G$, $m$, $T_\Phi$
		\State $T_\Psi=T_\Phi+T_G$, $M1= T_\Phi*diag(m)*T_\Psi'$	
		\State Initialize  $M1inv=M0inv$
		\State $n\gets e$, $m$, $T_\Phi$, $T_G$		
		\State {\bf if $n>0$}
		\State \ \ \  $P1=M0inv*(M0-M1)$
		\State \ \ \  Initialize $Q1=M0inv$
		\State \ \ \  {\bf for i:=1 to n}
		\State \ \ \  \ \ \  \ $Q1=P1*Q1$
		\State \ \ \  \ \ \  \ $M1inv=M1inv + Q1$
		\State \ \ \  {\bf end}
		\State {\bf end}
		\State $M_{m,\Phi,\Psi}^{-1} \gets M1inv$
	\end{algorithmic}
\end{algorithm}

\begin{algorithm} \small
	\caption{Iterative computation of $M_{m,\Phi,\Phi}^{-1}$ (M0inv) and $M_{m,\Phi,\Psi}^{-1}$ (M1inv):  [m,TPsi,M1,M2,M1inv,M2inv,n] = \textbf{Prop11InvMultOp}(c,r,TPhi,TPsi,m,e)	}
	\label{alg:Prop11}
	\begin{algorithmic}[1]	
	\State $T_\Psi\gets T_\Psi$, $T_\Phi$
	    \State $m\gets m$, $T_\Phi$, $T_\Psi$	    	
		\State Initialize  $M1inv=eye(r)$
		\State $n\gets e$, $m$, $T_\Phi$, $T_G$		
		\State {\bf if $n>0$}
		\State  \ \ \ $P1=eye(r)-M1$
		\State  \ \ \  Initialize $Q1=eye(r)$
		\State \ \ \  {\bf for i:=1 to n}
		\State \ \ \  \ \ \  \ $Q1=P1*Q1$
		\State \ \ \  \ \ \  \ $M1inv=M1inv + Q1$
		\State \ \ \  {\bf end}
		\State {\bf end}
		\State $M_{m,\Phi,\Psi}^{-1} \gets M1inv$
	\end{algorithmic}
\end{algorithm}

\vspace{.1in}

\newpage

{\bf Acknowledgements} 
The authors acknowledge support from the Austrian Science Fund (FWF) START-
project FLAME ('Frames and Linear Operators for Acoustical Modeling and
Parameter Estimation'; Y 551-N13). They thank to Z. Prusa for help with LTFAT.


\end{document}